\documentclass[12pt]{amsart}
\usepackage{amsfonts,latexsym,amsthm,amssymb}
\usepackage[all]{xy}

\input epsf

\newtheorem{theorem}{Theorem}[section]

{\theoremstyle{definition} \newtheorem{defin}[theorem]{Definition}}
{\theoremstyle{remark} \newtheorem{remark}[theorem]{Remark}
\newtheorem{example}[theorem]{Example}}

\newcommand{\Pbb}{{\mathbb{P}}}

\newcommand{\one}{1\hskip-3.5pt1}
\newcommand{\csm}{c_{\text{\rm SM}}}
\newcommand{\ssm}{s^\circ}

\title{Inclusion-exclusion and Segre classes}
\author{Paolo Aluffi}
\address{Max-Planck-Institut f\"ur Mathematik, Bonn, Germany}
\address{Florida State University, Tallahassee, Florida}

\makeindex

\begin{document}

\begin{abstract}
We propose a variation of the notion of Segre class, by forcing a
naive `inclusion-exclusion' principle to hold. The resulting class is
computationally tractable, and is closely related to
Chern-Schwartz-MacPherson classes. We deduce several general
properties of the new class from this relation, and obtain an
expression for the Milnor class of an arbitrary scheme in terms of
this class.
\end{abstract}

\maketitle


\section{Introduction}

Notwithstanding their fundamental r\^ole in modern intersection theory
(cf.~\cite{MR85k:14004}, Chapters~4 and~6), {\em Segre classes\/}
remain a somewhat esoteric concept: with a few notable exceptions,
they have been used more for foundational purposes than for actual
computations of concrete intersection products. This is due to the
effective inaccessibility of Segre classes: essentially no techniques
are known to compute the Segre class $s(Z,M)$ of a scheme $Z$ in a
scheme $M$, other than its raw definition; which is perhaps a little
too close to the ideal of $Z$ for comfort, in almost every
geometrically significant problem. The main virtue of the Segre
class---that is, its sensitivity to the fine structure of $Z$---turns
out being the main problem in handling it in concrete situations.

In this article we propose an a variation ($\ssm(Z,M)$,
Definition~\ref{ssmZ}) of the notion of Segre class of a subscheme $Z$
in a {\em nonsingular variety $M$,\/} by imposing on it an {\em
inclusion-exclusion principle,\/} which makes $\ssm(Z,M)$ well-behaved
with respect to naive set-theoretic operations. The class we define
does not work as a Segre class in a definition of an intersection
product, but shares with Segre classes several notable properties: our
$\ssm(Z,M)$ agrees with $s(Z,M)$ if $Z$ is nonsingular; and behaves
with respect to different embeddings of $Z$ in nonsingular varieties
or to smooth maps in precisely the same way the ordinary Segre class
would.

We collect these and other properties in Theorem~\ref{props}. We find
these observations very remarkable. First, it is very puzzling
that the definition of $\ssm(Z,M)$ makes sense at all. The definition
of $\ssm(Z,M)$ is given in terms of any collection of hypersurfaces
cutting out $Z$, and consists of a rather complicated combination of
conventional Segre classes. That the end-result should not depend on
the choices of the hypersurfaces must amount to massive cancellations
occurring among the Segre classes involved in the definition; to our
knowledge, there is no direct explanation for these cancellations. 
Second, the definition turns out to {\em only depend on the support
  of~$Z$;} given the sensitivity to scheme structure of
ordinary Segre classes, we find this fact rather astonishing, again
amounting to remarkable cancellations which must take into account and
then eliminate the contributions of nilpotents to the class. We
illustrate such `cancellations' with a couple of simple examples; the
reader is warmly invited to work out more complex ones.

Our perspective is that the properties listed in Theorem~\ref{props}
must reflect some unknown and powerful features of ordinary Segre
classes. We view Theorem~\ref{props} as `experimental evidence' for
these features, and the main purpose of this article is to advertise
this evidence. With this understood, it is remarkable that
Theorem~\ref{props} can be proved without uncovering or even being
able to state precisely these features. In fact, the proof of
Theorem~\ref{props} is essentially immediate once it is realized that
$\ssm(Z,M)$ is closely tied to another important class, defined for
all singular varieties. This relationship is exposed in
Theorem~\ref{ursprung}; Theorem~\ref{props} follows easily from this.

Unfortunately, being able to prove something is not the same as
understanding it. While we do prove Theorem~\ref{props}, the
properties of Segre classes which must be responsible
for it remain just as unknown after the fact. We believe that
establishing these properties would be exceedingly interesting. If
Segre classes were computable objects, then a great many problems in
enumerative geometry would become a routine exercise. Given the
(unreasonable?) relevance of enumerative geometry in recent
developments in algebraic geometry, clarifying the notion of Segre
class seems a very worthwhile goal.

This goal seems to us largely unmet---with the exception of the seminal
work of Steven Kleiman and his collaborators, in which Segre classes
play an important part (cf.~for example \cite{MR95c:14073} and
\cite{MR98g:14008}).

This is the first article in a series planned to explore
`inclusion-exclusion' phenomena in the theory of Segre classes. In
\cite{incexcII} we will propose a different variation on the theme of
Segre classes, also satisfying an inclusion-exclusion principle, and
yielding a simple computation of $\ssm(Z,M)$ in certain cases.\vskip 6pt

\noindent{\bf Acknowledgments} I thank the Max-Planck-Institut f\"ur
Mathematik in Bonn, Germany, for the hospitality and support.


\section{SM-Segre classes}\label{newseg}

Throughout this section, $M$ denotes a nonsingular variety (over
an algebraically closed field of characteristic~0, although this does
not seem to be essential). 

We will consider a proper subscheme $Z$ of $M$, and a finite family
$\{X_i\}_{i=1,\dots,r}$ of hypersurfaces cutting out $Z$ in $M$:
$$Z=X_1\cap\cdots\cap X_r\quad;$$
note that we are putting no restriction on $r$, and no other
restrictions on the hypersurfaces (they may be nonreduced, there may
be repetitions in the list, etc.). In fact, the requirement on the
hypersurfaces will be further relaxed later on, cf.~Remark~\ref{supsuf}.

For a hypersurface $X$, we define its {\em SM-Segre class\/} as
follows\footnote{{\em SM} is supposed both to evoke the connection
  with Schwartz-MacPherson classes, which is the key to the main
  properties of the class, and the fact that we view this notion
  as a `{\em sm\/}ooting' of the notion of conventional Segre
  class}. Let $Y$ be the singularity subscheme of $X$, that is, the
subscheme locally defined by the partial derivatives of a local
generator for the ideal of $X$. 

\begin{defin}\label{ssmX}
The SM-Segre class of $X$ in $M$ is the class
$$\ssm(X,M)=s(X,M)+c(\mathcal O(X))^{-1}\cap(s(Y,M)^\vee
\otimes_M \mathcal O(X))\quad.$$
\end{defin}

This definition uses notations---e.g., for the tensor of a rational
equivalence class by a line bundle---introduced in \cite{MR96d:14004},
Def.~2; in more conventional (but less manageable) terms, the component
of dimension $m$ in $\ssm(X,M)$ is
$$s(X,M)_m+(-1)^{n-m}\sum_{j=0}^{n-m}\binom {n-m}j X^j\cdot
s(Y,M)_{m+j}$$
where $n=\dim M$; a formula that is reminiscent of classical residual
intersection formulas (cf.~\cite{MR85k:14004}, Prop.~9.2). In
\cite{MR96d:14004} and \cite{MR2001i:14009} the class $\ssm(X,M)$ is
denoted $s(X\setminus Y,M)$; no analogs for schemes other than
hypersurfaces were considered there.

As given in Definition~\ref{ssmX}, the SM-Segre class {\em of a
hypersurface $X$\/} lives naturally in the Chow group $A_*X$ of
$X$. In view of the upgrade to arbitrary subschemes $Z$ of $M$ that
follows, however, it will be more natural to view it for the moment as
an element of the Chow group $A_*M$ of the ambient nonsingular
variety. We will omit such evident push-forwards from our notations,
as well as evident pull-backs.

\begin{defin}\label{ssmZ}
Let $Z$ be a proper subscheme of $M$, and let $X_1,\dots,X_r$ be
hypersurfaces cutting out $Z$:
\begin{equation}\label{equal}
Z=X_1\cap\cdots\cap X_r\quad.
\end{equation}
Then the {\em SM-Segre class of $Z$ in $M$\/} is obtained by applying
inclusion-exclusion to the SM-Segre classes of the hypersurfaces
$X_i$. Explicitly, we set
$$\ssm(Z,M)=\sum_{s=1}^r (-1)^{s-1} \sum_{i_1<\dots<i_s}
\ssm(X_{i_1}\cup\cdots\cup X_{i_s},M)\quad\in A_*M\quad.$$
\end{defin}
Here $X_{i_1}\cup\cdots\cup X_{i_s}$ is the hypersurface whose ideal
is the product of the ideals of $X_{i_1},\dots,X_{i_s}$ (so $X\cup
X\ne X$!); but see Remark~\ref{supsuf} below. Also, while we are
defining the class in $A_*M$ at this stage, Theorem~\ref{ursprung}
will imply that it is the image of class naturally defined on $Z$
itself; also cf.~Remark~\ref{indlim}.

The terminology {\em inclusion-exclusion\/} is adapted from the
similar-looking formula for the number of elements in the
intersection $Z$ of $r$ finite sets $X_1,\dots,X_r$:
$$\# Z=\sum_{s=1}^r (-1)^{s-1} \sum_{i_1<\dots<i_s} \#
(X_{i_1}\cup\cdots\cup X_{i_s})\quad;$$
this is immediately proved by induction on $r$.

Concerning Definition~\ref{ssmZ}, the reader should now expect a proof
that $\ssm(Z,M)$ does not depend on the specific choice of
hypersurfaces cutting out $Z$. We isolate this and other properties of
the definition in the following statement.

\begin{theorem}\label{props}
\begin{enumerate}
\item\label{choices} The definition of $\ssm(Z,M)$ is independent of
  the choices.
\item\label{support} If $Z_{\text red}$ is the support of $Z$, then
  $\ssm(Z,M)=\ssm(Z_{\text red},M)$.
\item\label{ns} If $Z$ is nonsingular, then
  $\ssm(Z,M)=s(Z,M)=c(N_ZM)^{-1}\cap [Z]$.
\item\label{closed} If $Z\subset M\subset M'$, where $M'$ is a
  nonsingular variety and $M\subset M'$ is a closed embedding, then
$$\ssm(Z,M')=c(N_MM')^{-1}\cap \ssm(Z,M)\quad.$$
\item\label{open} If $Z\subset U\subset M$, where $U\stackrel j\to M$
  is an open embedding, then $\ssm(Z,U)=j^* \ssm(Z,M)$.
\item\label{smooth} More generally: if $p: M'\to M$ is a smooth
  morphism, $Z\subset M$   a subscheme, and $Z'=p^{-1}(Z)$ its inverse
  image, then $\ssm(Z',M')=p^* \ssm(Z,M)$.
\item\label{full} The class $\ssm(Z,M)$ satisfies a full
  inclusion-exclusion principle, in the sense that if $Z_1,\dots, Z_r$
  are subschemes of $M$ such that $Z=Z_1\cap\cdots\cap Z_r$, then
$$\ssm(Z,M)=\sum_{s=1}^r (-1)^{s-1} \sum_{i_1<\dots<i_s}
\ssm(Z_{i_1}\cup\cdots\cup Z_{i_s},M)\quad,$$
where $Z_{i_1}\cup\cdots\cup Z_{i_s}$ is the subscheme of $M$ whose
ideal is the product of the ideals of $Z_{i_1},\dots, Z_{i_s}$.
\end{enumerate}
\end{theorem}

We are separating these statements from their {\em Ursprung,\/} which
is Theorem~\ref{ursprung} below, in the attempt to highlight them
independently of our technical bias. Theorem~\ref{props} will follow from
Theorem~\ref{ursprung}, but we strongly feel that these statements are
of substantial independent interest, and call for a straightforward
intersection-theoretic proof; with the exclusion of
parts~\ref{open}.-\ref{full}., which are formal exercises (left to the
reader), we do not know such a proof.

We delay the (rather anticlimatic) proof of Theorem~\ref{props} until
the next section. The rest of this section is taken by several
comments meant to further highlight the content of the theorem.

\begin{remark}\label{supsuf}
By part~\ref{support}., the equality (\ref{equal}) in
Definition~\ref{ssmZ} need only hold set-theo\-retically. By the
same token, we can in fact replace the unions $X_{i_1}\cup\cdots\cup
X_{i_s}$ in that definition, or the union $Z_{i_1}\cup\cdots\cup
Z_{i_s}$ in part \ref{full}., by any other schemes supported on such
unions (for example, we could take the ideals defined by the {\em
  intersection\/} of the ideals, rather than their product, or take
radicals throughout).
\qed\end{remark}

\begin{remark}\label{consns}
In view of parts~\ref{ns}.~and \ref{closed}, it is consistent to define
$\ssm(M,M)=[M]$.
\qed\end{remark}

\begin{remark}\label{indlim}
The formulas in part~\ref{closed}.~and~\ref{open}., ruling the
behavior of the SM-Segre class under different embeddings of $Z$, hold
in precisely the same terms for the ordinary Segre class (this follows
from \cite{MR85k:14004}, Example~4.2.6(a)). In fact, part~\ref{closed}
suggests that $\ssm(Z,M)$ should be, up to a correction term $c(TM)$,
a class intrinsic of $Z$ (and living in $A_*Z$). This is precisely the
case, as will follow from Theorem~\ref{ursprung}.
\qed\end{remark}

\begin{remark}\label{local}
By part~\ref{ns}.~and \ref{open}., if $Z$ is reduced then the
difference between $\ssm(Z,M)$ and $s(Z,M)$ is supported within the
singular locus $Z_s$ of $Z$. Indeed, letting $U=M\setminus Z_s$, the
classes $\ssm(Z,M)$ and $s(Z,M)$ restrict (by part~\ref{open}.) to the
same (by part~\ref{ns}.) class in $A_*U$, so the difference comes from
$A_*Z_s$ by \cite{MR85k:14004}, Proposition~1.8.
\qed\end{remark}

\begin{remark}\label{birinv}
The behavior of $\ssm$ under smooth morphisms, prescribed by
part~\ref{smooth}., also matches the behavior of conventional Segre
classes. In fact, for Segre classes the same behavior extends to the
more general case of {\em flat\/} maps, by \cite{MR85k:14004},
Proposition~4.2(b); this is not so for SM-Segre classes. Another key
property of Segre classes, that is, birational invariance
(cf.~Proposition~4.2(a) in \cite{MR85k:14004}) also does not hold for
SM-Segre classes.

On the other hand, by part~\ref{support}.~of Theorem~\ref{props}, the
stated equality for SM-Segre classes under smooth morphisms holds as
soon as $Z'$ equals $p^{-1}(Z)$ set-theoretically; for the ordinary
Segre class, this equality has to hold scheme-theoretically.
\qed\end{remark}

To give a sense of what might go into a `direct argument' as
envisioned above, here is a direct proof of a very particular case of
the innocuous-looking part~\ref{ns}.

\begin{proof} We will prove part~\ref{ns}.~in the particular case
in which $Z$ is the complete intersection of two transversal
nonsingular hypersurfaces $X_1$, $X_2$. In this case,
Definition~\ref{ssmX} gives
$$\ssm(X_1,M)=\frac{[X_1]}{1+X_1}\quad,\quad
\ssm(X_2,M)=\frac{[X_2]}{1+X_2}$$
(where we employ the convenient shorthand
$\frac{[X]}{1+X}=[X]-[X^2]+[X^3]-\dots =s(X,M)$). As for $X_1\cup
X_2$, it is easily verified that $Z$ is itself the singularity
subscheme of this hypersurface; hence
$$\aligned
\ssm(X_1\cup X_2,M) &=\frac{[X_1+X_2]}{1+X_1+X_2} + \frac 1{1+X_1+X_2}
\left(s(Z,M)^\vee\otimes \mathcal O(X_1+X_2)\right)\\
&=\frac{[X_1+X_2]}{1+X_1+X_2} + \frac 1{1+X_1+X_2}
\left(\frac{[X_1 X_2]}{(1-X_1)(1-X_2)}\otimes \mathcal O(X_1+X_2)\right)\\
&=\frac{[X_1+X_2]}{1+X_1+X_2} + \frac 1{1+X_1+X_2}
\frac{[X_1 X_2]}{(1+X_2)(1+X_1)}
\endaligned$$
where we have used Proposition~1 in \cite{MR96d:14004}. Thus,
according to Definition~\ref{ssmZ}:
$$\aligned
\ssm(Z,M)&=\frac{[X_1]}{1+X_1}+\frac{[X_2]}{1+X_2}
-\left(\frac{[X_1+X_2]}{1+X_1+X_2} + \frac 1{1+X_1+X_2} \frac{[X_1
    X_2]}{(1+X_2)(1+X_1)}\right)\\
&= \frac{[X_1 X_2]}{(1+X_1)(1+X_2)}
\endaligned$$
by trivial formal manipulations. Since $Z$ is the complete intersection
of $X_1$ and $X_2$, the right-hand-side is $s(Z,M)$, as prescribed by
part~\ref{ns} of Theorem~\ref{props}.
\end{proof}

The whole of Theorem~\ref{props}, especially the crucial
parts~\ref{choices}.~and \ref{support}., ought to have a similarly
direct explanation, but none is available to us at this time. For
Theorem~\ref{props} to hold, drastic cancellations (of which the
simplifications occurring in the proof presented a moment ago must be
the simplest instance) must be at work behind the scenes. This can be
also be observed in any concrete computation of SM-Segre classes; we
give two simple examples for illustration purposes.

\begin{example}\label{xxx} The curve $C$ consisting of two lines
  meeting at a point in $\Pbb^3$ can be realized as the intersection
  of a nonsingular quadric $Q$ and a tangent plane $H$. Applying the
  definition of SM-Segre classes (noting that $Y=\emptyset$ for a
  nonsingular hypersurface) gives
$$\ssm(H,\Pbb^3)=s(H,\Pbb^3)=[\Pbb^2]-[\Pbb^1]+[\Pbb^0]$$
$$\ssm(Q,\Pbb^3)=s(Q,\Pbb^3)=2[\Pbb^2]-4[\Pbb^1]+8[\Pbb^0]\quad.$$
The singularity subscheme $Y$ of $H\cap Q$ consists of the two lines, with
an embedded point at the intersection point. A blow-up computation gives 
$$s(Y,\Pbb^3)=2[\Pbb^1]-4[\Pbb^0]$$
yielding, according to Definition~\ref{ssmX},
$$\ssm(H\cup Q,\Pbb^3)=3[\Pbb]^2-7[\Pbb]^1+13[\Pbb^0]\quad.$$
Thus (Definition~\ref{ssmZ}):
$$\aligned
\ssm(C,\Pbb^3)&=([\Pbb^2]-[\Pbb^1]+[\Pbb^0])+(2[\Pbb^2]-4[\Pbb^1]
+8[\Pbb^0])-(3[\Pbb]^2-7[\Pbb]^1+13[\Pbb^0])\\
&=2[\Pbb^1]-4[\Pbb^1]\quad.
\endaligned$$
On the other hand, $C$ is itself a hypersurface, in $M'=\Pbb^2$, with
singularity subscheme a point. By Definition~\ref{ssmX},
$$\ssm(C,\Pbb^2)=\frac{[C]}{1+C}+[\Pbb^0]=2[\Pbb^1]-3[\Pbb^0]\quad.$$
As prescribed by part \ref{closed}.~of Theorem~\ref{props},
$$\ssm(C,\Pbb^3)=c(N_{\Pbb^2}\Pbb^3)^{-1}\cap \ssm(C,\Pbb^2)\quad.$$
\end{example}

\begin{example}\label{embed}
Now let $C$ be the subscheme defined by the homogeneous ideal
$(xy,x^2)$ in $\Pbb^2$; that is, a line with an embedded point. Thus,
$C$ is the intersection of $X_1$ and $X_2$, where $X_1$ is a union of
two lines, and $X_2$ is the double line with ideal $(x^2)$. By
Example~\ref{xxx} we have
$$\ssm(X_1,\Pbb^2)=2[\Pbb^1]-3[\Pbb^0]\quad.$$
The singularity subscheme of $X_2$ consist of a (single) line, so
Definition~\ref{ssmX} gives
$$\aligned
\ssm(X_2,\Pbb^2)&=\frac{[X_2]}{1+X_2}+\frac
1{1+X_2}\left((s(\Pbb^1,\Pbb^2)^\vee
\otimes_{\Pbb^2}\mathcal(X_2)\right)\\
&=2[\Pbb^1]-4[\Pbb^0]-([\Pbb^1]-3[\Pbb^0])=[\Pbb^1]-[\Pbb^0]\quad.
\endaligned$$
Next, $X_1\cup X_2$ has ideal $(x^3y)$, hence its singularity
subscheme has ideal $(x^2y,x^3)$. Computing (conventional) Segre
classes and using again Definition~\ref{ssmX} gives
$$\ssm(X_1\cup X_2,\Pbb^2)=2[\Pbb^1]-3[\Pbb^0]\quad.$$
Thus
$$\ssm(C,\Pbb^2)=(2[\Pbb^1]-3[\Pbb^0])+([\Pbb^1]-[\Pbb^0])-(2[\Pbb^1]
-3[\Pbb^0])=[\Pbb^1]-[\Pbb^0]\quad.$$
This is in agreement with parts \ref{support}.~and \ref{ns}.~of
Theorem~\ref{props}: the support of $C$ is $\Pbb^1$, hence
$$\ssm(C,\Pbb^2)=\ssm(\Pbb^1,\Pbb^2)=s(\Pbb^1,\Pbb^2)
=[\Pbb^1]-[\Pbb^0]$$
as obtained `by hand'.
\end{example}


\section{The proof, and comments about Milnor classes}\label{miln}

Theorem~\ref{props} can be proved rather indirectly, by applying one
previous result and a piece of a well-established theory. Once more,
our goal in this article is really to suggest that Theorem~\ref{props}
ought to have a direct proof in terms of the theory of Segre classes,
and that finding this argument would tell us something very
interesting about Segre classes. Our proof does not shed much light in
this direction.

We denote by $\csm$ the {\em Chern-Schwartz-MacPherson\/} class,
cf.~\cite{MR50:13587} (and \cite{MR91h:14010} for a treatment over an
arbitrary algebraically closed field of characteristic~0).

\begin{theorem}\label{ursprung}
Let $Z$ be a subscheme of a nonsingular variety $M$. Then
$$\csm(Z_{\text{red}})=c(TM)\cap \ssm(Z,M)\quad.$$
\end{theorem}

\begin{proof}
The Chern-Schwartz-MacPherson class satisfies inclusion-exclusion.
Indeed, if 
$$Z=X_1\cup\cdots\cup X_n\quad,$$
then one easily checks that
$$\one_Z=\sum_{s=1}^r (-1)^{s-1} \sum_{i_1<\dots<i_s} \one_{X_{i_1}
  \cup\cdots\cup X_{i_s}}\quad,$$
where $\one_X$ denotes the constructible function which is 1 at points
  of $X$, and 0 outside of $X$; and note $\one_X=\one_{X_{\text{red}}}$
  trivially. Applying MacPherson's natural transformation gives then
$$\csm(Z_{\text{red}})=\sum_{s=1}^r (-1)^{s-1} \sum_{i_1<\dots<i_s}
  \csm((X_{i_1}\cup\cdots\cup X_{i_s})_{\text{red}})$$
in $A_*M$. This observation reduces the statement of the theorem to
  the case in which $Z$ is a hypersurface in a nonsingular variety,
  which is shown in \cite{MR2001i:14009} (Theorem~I.1 and Corollary~II.2).
\end{proof}

Theorem~\ref{props} follows immediately as a corollary; the details
are left to the reader.

We also note that since $\ssm(Z,M)=c(TM)^{-1}\cap
\csm(Z_{\text{red}})$, it follows that $\ssm(Z,M)$ is naturally the
image of a class in $A_*Z$, although this is far from clear from
Definition~\ref{ssmZ} (cf.~Remark~\ref{indlim}).\vskip 6pt

We end with comments regarding the so-called {\em Milnor class\/} of a
scheme.

Theorem~\ref{ursprung} writes the Chern-Schwartz-MacPherson class in
a way that the informed reader will recognize immediately as a
companion to the definition of an intrinsic class given by William
Fulton (cf.~\cite{MR85k:14004}, 4.2.6 for the definition of this class
and of a kindred notion, the {\em Fulton-Johnson\/} class): if $Z$ is
a subscheme of a nonsingular variety $M$, Fulton's class is defined by 
$$c_F(Z)=c(TM)\cap s(Z,M)\quad.$$

According to several authors, but with slightly different conventions
of sign and context, the {\em Milnor class\/} measures the difference
between the Chern-Schwartz-MacPher\-son class of a variety and other
classes such as Fulton's or Fulton-Johnson's. Determining this
discrepancy has been identified as a Verdier-Riemann-Roch type
problem, cf.~\cite{MR2000e:14011}. The hypersurface case is rather
well understood, cf.~\cite{MR1795550}; complete intersections are
treated in \cite{MR2001g:32062} and in the recent
\cite{math.AG/0202175}.

Our viewpoint on the problem is perhaps a little different from the
one taken by these authors. For us, the Milnor class of $Z$ should be
measured by a Segre-class type of invariant defined on the singularity
subscheme of $Z$; in fact our motivation in pursuing any such formula
is precisely to learn something new about Segre classes. In this
sense, the problem seems wide open for anything but hypersurfaces.

In the context of the present article, the relevant remark is the
following consequence of Theorem~\ref{ursprung}:
the Milnor class of a reduced scheme $Z$ embedded in a nonsingular
variety $M$ is (up to sign)
$$c(TM)\cap m(Z,M)\quad,$$
where 
$$m(Z,M):=\ssm(Z,M)-s(Z,M)\quad.$$
Note that the class $m(Z,M)$ is localized on the singularities of $Z$
(see Remark~\ref{local}). 

From our perspective, the main problem in the study of Milnor classes
is the explicit determination of $m(Z,M)$ in terms of conventional
intersection-theoretic operations.

\begin{example} 
If $X$ is a hypersurface, with singularity subscheme $Y$, then
$$m(X,M)=c(\mathcal O(X))^{-1}\cap(s(Y,M)^\vee \otimes_M \mathcal
O(X))\quad.$$
\end{example}
This formula, now incorporated in the definition of $\ssm$, is
the main result of \cite{MR96d:14004} and \cite{MR2001i:14009}.
Similarly explicit formulas for $m(Z,M)$ for more general $Z$, in
terms of Chern/Segre-class type invariants of the singularity
subscheme of $Z$, would be highly desirable. To our knowledge, it is
not even known whether $m(Z,M)$---and hence the Milnor class of
$Z$---is determined by the singularity subscheme of $Z$; {\em even
  when $Z$ is a complete intersection of codimension~2.}




\begin{thebibliography}{BLSS99}

\bibitem[Alu94]{MR96d:14004}
Paolo Aluffi.
\newblock Mac{P}herson's and {F}ulton's {C}hern classes of hypersurfaces.
\newblock {\em Internat. Math. Res. Notices}, (11):455--465, 1994.

\bibitem[Alu99]{MR2001i:14009}
Paolo Aluffi.
\newblock Chern classes for singular hypersurfaces.
\newblock {\em Trans. Amer. Math. Soc.}, 351(10):3989--4026, 1999.

\bibitem[Alu02]{incexcII}
Paolo Aluffi.
\newblock {Inclusion-Exclusion and Segre classes, II}.
\newblock In preparation, 2002.

\bibitem[BLSS99]{MR2001g:32062}
Jean-Paul Brasselet, Daniel Lehmann, Jos{\'e} Seade, and Tatsuo Suwa.
\newblock Milnor numbers and classes of local complete intersections.
\newblock {\em Proc. Japan Acad. Ser. A Math. Sci.}, 75(10):179--183, 1999.

\bibitem[Ful84]{MR85k:14004}
William Fulton.
\newblock {\em Intersection theory}.
\newblock Springer-Verlag, Berlin, 1984.

\bibitem[Ken90]{MR91h:14010}
Gary Kennedy.
\newblock Mac{P}herson's {C}hern classes of singular algebraic varieties.
\newblock {\em Comm. Algebra}, 18(9):2821--2839, 1990.

\bibitem[Kle94]{MR95c:14073}
Steven~L. Kleiman.
\newblock A generalized {T}eissier-{P}l\"ucker formula.
\newblock In {\em Classification of algebraic varieties (L'Aquila, 1992)},
  pages 249--260. Amer. Math. Soc., Providence, RI, 1994.

\bibitem[KT96]{MR98g:14008}
Steven Kleiman and Anders Thorup.
\newblock Mixed {B}uchsbaum-{R}im multiplicities.
\newblock {\em Amer. J. Math.}, 118(3):529--569, 1996.

\bibitem[Mac74]{MR50:13587}
R.~D. MacPherson.
\newblock Chern classes for singular algebraic varieties.
\newblock {\em Ann. of Math. (2)}, 100:423--432, 1974.

\bibitem[PP01]{MR1795550}
Adam Parusi{\'n}ski and Piotr Pragacz.
\newblock Characteristic classes of hypersurfaces and characteristic cycles.
\newblock {\em J. Algebraic Geom.}, 10(1):63--79, 2001.

\bibitem[Sch]{math.AG/0202175}
Joerg Sch\"urmann.
\newblock {A generalized Verdier-type Riemann-Roch theorem for
  Chern-Schwartz-MacPherson classes}, arXiv:math.AG/0202175.

\bibitem[Yok99]{MR2000e:14011}
Shoji Yokura.
\newblock On a {V}erdier-type {R}iemann-{R}och for
  {C}hern-{S}chwartz-{M}ac{P}herson class.
\newblock {\em Topology Appl.}, 94(1-3):315--327, 1999.
\newblock Special issue in memory of B. J. Ball.

\end{thebibliography}
\end{document}